\documentclass[12pt]{article}
\usepackage[T1]{fontenc}
\usepackage{amsmath}
\usepackage{amssymb}
\usepackage{amsthm}
\usepackage{pst-all}
\newtheorem{theorem}{Theorem}[section]
\newtheorem{definition}[theorem]{Definition}
\newtheorem{lemma}[theorem]{Lemma}
\newtheorem{proposition}[theorem]{Proposition}
\newtheorem{corollary}[theorem]{Corollary}

\def\fl#1{\smash{\mathop{\hbox to 12mm{ \rightarrowfill\ }}\limits^{\textstyle #1}}}
\newcommand{\qtt}{\Lambda}
\newcommand{\fin}{\hfill$\diamond$ \\}
\newcommand{\lhp}{\textrm{{\huge (}}}
\newcommand{\lbp}{\textrm{{\Large (}}}
\newcommand{\rhp}{\textrm{{\huge )}}}
\newcommand{\rbp}{\textrm{{\Large )}}}
\newcommand{\codim}{\textrm{codim}}

\title{On Alexander modules and Blanchfield forms of null-homologous knots in rational homology spheres}
\author{Delphine Moussard}
\date{ }

\begin{document}

\maketitle

\begin{abstract}
In this article, we give a classification of Alexander modules of null-homologous knots in rational homology spheres. 
We characterize these modules $\mathcal{A}$ equipped with their Blanchfield forms $\phi$, and the modules $\mathcal{A}$ such that
there is a unique isomorphism class of $(\mathcal{A},\phi)$, and we prove that for the other modules $\mathcal{A}$, 
there are infinitely many such classes. We realise all these $(\mathcal{A},\phi)$ by explicit knots in $\mathbb{Q}$-spheres.

\ \\
MSC : 57M25 57M27 57N10 57N65

\ \\
Keywords : Homology sphere; Surgery; Alexander module; Equivariant linking number; Equivariant linking matrix; Blanchfield pairing;
Blanchfield form. 
\end{abstract}
\newpage
\tableofcontents

\section{Introduction}

  \subsection{Short introduction}

The $\mathbb{Q}[t,t^{-1}]$-modules appearing as Alexander modules of knots in $S^3$ were determined by Levine \cite{Lev}.
They are all the $ \bigoplus_{i=1}^k \mathbb{Q} [t,t^{-1}]/(P_i)$, where the 
$P_i$ are symmetric polynomials ($P_i(t)=P_i(t^{-1})$) in $\mathbb{Z} [t,t^{-1}]$ such that $P_i(1)=1$.
We will generalize this classification to null-homologous knots in rational homology spheres.
The Alexander modules of these knots are finitely generated $\mathbb{Q}[t,t^{-1}]$-torsion-modules on which 
$x\mapsto (1-t)x$ defines an isomorphism.
The equivariant linking number defined in the infinite cyclic covering induces a hermitian form
on these Alexander modules, called the Blanchfield form since Blanchfield showed that it is non degenerate \cite{Blanch}.
We will show that these properties characterize the Alexander modules $\mathcal{A}$ equipped with their Blanchfield forms $\phi$.
We characterize these modules $\mathcal{A}$ equipped with their Blanchfield forms $\phi$, and the modules $\mathcal{A}$ such that
there is a unique isomorphism class of $(\mathcal{A},\phi)$, and we prove that for the other modules $\mathcal{A}$, 
there are infinitely many such classes.

In \cite{Kricker}, Kricker defined a rational lift of the Kontsevich integral for null-homologous knots in $\mathbb{Q}$-spheres. 
Garoufalidis and Rozansky \cite{G-R} introduced a filtration of the space of these knots by null moves.
Null moves preserve the isomorphism classes of $(\mathcal{A},\phi)$. Garoufalidis and Rozansky computed the graded space 
of the filtration for the isomorphism class of the unknot. 
This computation allowed Lescop to prove that the invariant she constructed in \cite{Lescop} is equivalent 
to the 2-loop part of the Kricker lift for knots with trivial Alexander modules. The work contained here should be useful 
to generalize this result.

  \subsection{Statement of the results}

Except otherwise mentioned, all manifolds will be compact and oriented, 
and all manifolds of dimension 2 or 3 will be connected. We set $\Lambda=\mathbb{Q}[t,t^{-1}]$.

A \emph{rational homology 3-sphere}, or \emph{$\mathbb{Q}$-sphere}, is a 3-manifold, 
without boundary, which has the same rational homology as the standard sphere $S^3$.
In such a $\mathbb{Q}$-sphere $M$, a \emph{null-homologous knot} $K$ is a knot whose 
homology class in $H_1(M;\mathbb{Z})$ is zero. 
Let $T(K)$ be a tubular neighborhood of $K$. The \emph{exterior} of $K$ is 
$X=M\setminus Int(T(K))$. Consider the projection $\pi : \pi_1(X) \to \frac{H_1(X;\mathbb{Z})}{torsion} \cong \mathbb{Z}$,
and the covering map $p : \tilde{X} \to X$ associated with its kernel. The covering $\tilde{X}$ is the \emph{infinite cyclic covering} 
of $X$. The automorphism group of the covering, $Aut(\tilde{X})$, is isomorphic to $\mathbb{Z}$. It acts on 
$H_1(\tilde{X};\mathbb{Q})$. Denoting the action of a generator $\tau$ of $Aut(\tilde{X})$ as the multiplication by $t$, 
we get a structure of $\qtt$-module on 
$H_1(\tilde{X};\mathbb{Q})$. This $\qtt$-module is called the \emph{Alexander module} of $K$, and we will 
denote it by $\mathcal{A}(K)$.

On the Alexander module $\mathcal{A}(K)$, one can define the \emph{Blanchfield form}, or \emph{equivariant linking pairing},  
$\phi_K : \mathcal{A}(K)\times\mathcal{A}(K) \to \frac{\mathbb{Q}(t)}{\qtt}$, as follows. First define the equivariant linking number of two knots.
\begin{definition}
 Let $J_1$ and $J_2$ be two links in $\tilde{X}$ such that $J_1\cap \tau^k(J_2)=\emptyset$ for all $k\in\mathbb{Z}$.
 Let $\delta(t)$ be the annihilator of $\mathcal{A}(K)$.
 Then $\delta(\tau)J_1$ and $\delta(\tau)J_2$ are rationally null-homologous knots. The \emph{equivariant linking number} of $J_1$ and $J_2$ is 
 $$lk_e(J_1,J_2)=\frac{1}{\delta(t)\delta(t^{-1})}\sum_{k\in\mathbb{Z}}lk(\delta(\tau)J_1,\tau^k(\delta(\tau)J_2))t^k.$$
\end{definition}
One can easily see that $lk_e(J_1,J_2)\in\frac{1}{\delta(t)}\qtt$, and $lk_e(J_2,J_1)(t)=lk_e(J_1,J_2)(t^{-1})$.
Now, if $\gamma$ (resp. $\eta$) is the homology class of $J_1$ (resp. $J_2$) in $\mathcal{A}(K)$, define $\phi_K(\gamma,\eta)$ by :
$$\phi_K(\gamma,\eta)=lk_e(J_1,J_2)\ mod\ \qtt$$
The Blanchfield form is hermitian ($\phi(\gamma,\eta)(t)=\phi(\eta,\gamma)(t^{-1})$ for all $\gamma,\eta\in\mathcal{A}(K)$),
and non degenerate : $\phi_K(\gamma,\eta)=0$ for all $\eta\in\mathcal{A}(K)$ implies $\gamma=0$.

\begin{proposition} \label{prop1}
 Let $(\mathcal{A}(K),\phi_K)$ be the Alexander module and the Blanchfield form of a null-homologous knot $K$ in a $\mathbb{Q}$-sphere $M$.
  \begin{enumerate}
   \item The module $\mathcal{A}(K)$ is a finitely generated $\qtt$-torsion-module.
   \item The map $x\mapsto(1-t)x$ defines an isomorphism of $\mathcal{A}(K)$.
   \item The form $\phi_K$ is hermitian and non degenerate.
  \end{enumerate}
\end{proposition}

Blanchfield showed (3.), we will show (1. 2.) in Section \ref{sechpm}.
We will also show the following description of $(\mathcal{A}(K),\phi_K)$.
\begin{theorem} \label{thdecomposition}
 If $(\mathcal{A},\phi)$ satisfy Conditions (1. 2. 3.), then $\mathcal{A}$ is a direct sum, orthogonal with respect to $\phi$,
 of submodules of these two kinds :
\begin{itemize}
 \item $\frac{\qtt}{(\pi^n)}\gamma$, with $\pi$ prime and symmetric, or $\pi=t+2+t^{-1}$, $n>0$,
  and $\phi(\gamma,\gamma)=\frac{P}{\pi^n}$, $P$ symmetric and prime to $\pi$.
 \item $\frac{\qtt}{(\pi^{n})}\gamma_1 \oplus \frac{\qtt}{(\overline{\pi}^{n})}\gamma_2$,
  with either $\pi$ prime, non symmetric, $\pi(-1)\neq 0$, $n>0$, or $\pi=1+t$, $n$ odd, and in both cases 
 $\phi(\gamma_1,\gamma_2)=\frac{1}{\pi^n}$, $\phi(\gamma_i,\gamma_i)=0$ for $i=1,2$.
\end{itemize}
\end{theorem}

This description will allow us to show the reciprocal result of Proposition \ref{prop1}.
\begin{theorem} \label{characterization}
 If $(\mathcal{A},\phi)$ satisfy Conditions (1. 2. 3.), then there is a null-homologous knot $K$ in a $\mathbb{Q}$-sphere $M$
 such that $(\mathcal{A}(K),\phi_K)$ is isomorphic to $(\mathcal{A},\phi)$.
\end{theorem}
We will give an explicit construction of a knot for the two kinds of modules of Theorem \ref{thdecomposition}, and we will
get the general case by using connected sums.

Considering only the Alexander modules, we have the following classification.
\begin{theorem} \label{Amodules}
Given a family $(\delta_1, \dots, \delta_p)$ of polynomials in $\qtt$ such that
$\delta_{i+1}|\delta_i$ for $1\leq i<p$, the module $\bigoplus_{i=1}^p \frac{\qtt}{\delta_i(t)}$
is the Alexander module of a null-homologous knot in a $\mathbb{Q}$-sphere if and only if 
the $\delta_i$ satisfy the following conditions :
\begin{itemize}
	\item $\delta_i(1)\neq 0$ for $1\leq i\leq p$,
	\item $\delta_i(t^{-1})= t^{q_i}\delta_i(t)$, with $q_i\in\mathbb{Z}$, for $1\leq i\leq p$,
	\item if, for $1\leq i\leq p$, $m_i$ denotes the multiplicity of $-1$ as a root of $\delta_i$,
		then, for any odd integer $m$, the number of indices $i$ such that $m_i=m$ is even.
\end{itemize}
\end{theorem}
Then the \emph{Alexander polynomial} of $K$ is the order of $\mathcal{A}(K)$, $\Delta=\amalg_{i=1}^p \delta_i$.
The \emph{annihilator} of $\mathcal{A}(K)$ is $\delta=\delta_1$.
This result can be viewed as a corollary of the previous ones, but we will show it first, independently of the existence of the 
Blanchfield form.

Moreover, using algebraic number theory, we will show :
\begin{theorem} \label{thBf}
 Let $\mathcal{A}=\bigoplus_{i=1}^p \frac{\qtt}{\delta_i(t)}$, with $\delta_{i+1}|\delta_i$ for $1\leq i<p$.
 Let $\Delta=\amalg_{i=1}^p \delta_i$ be the order of $\mathcal{A}$. For $1\leq i\leq p$, 
 let $m_i$ denote the multiplicity of $-1$ as a root of $\delta_i$. If $\Delta$ has no prime and symmetric divisor, 
 and if all $m_i$ are odd, then there is a unique isomorphism class of non degenerate hermitian forms
 $\phi:\mathcal{A}\times\mathcal{A}\to\mathbb{Q}(t)/\qtt$. Otherwise,
 there are infinitely many isomorphism classes of such forms.
\end{theorem}

  \subsection{Plan of the article}

The first section is devoted to the proof of Theorem \ref{Amodules}. We introduce surgery presentations of knots in Section 2.1, 
and show in Section 2.2 that the associated equivariant linking matrices are presentation matrices of the Alexander module. 
We deduce the properties of the Alexander modules in Section 2.3. In Section 2.4, we realize any hermitian matrix with 
coefficients in $\qtt$ whose determinant does not vanish at $t=1$ as an equivariant linking matrix. We conclude 
the proof of Theorem \ref{Amodules} in Section 2.5.

Section 3.1 gives an expression of the Blanchfield form in terms of an equivariant linking matrix 
associated with a surgery presentation of the knot. 
Then, we prove Theorem \ref{thdecomposition} in Section 3.2, Theorem \ref{characterization} in Section 3.3, 
and Theorem \ref{thBf} in Section 3.4.

  \subsection{Acknowledgements}

I wish to thank Greg Kuperberg and Ga$\ddot{\textrm{e}}$l R\'emond for their help and ideas concerning the use of algebraic 
number theory in the last section.
I am very thankful to my advisor, Christine Lescop, for her precious guidance and advice.

\section{Classification of Alexander modules}

  \subsection{Surgery presentation of a knot}

Let $M_0$ be a 3-manifold such that $H_1(M_0;\mathbb{Q})=0$. Consider 
a link $L=\amalg_{i\in I} J_i$ in $M_0$, with knot components $J_i$. For each $i\in I$, let $T(J_i)$ be 
a tubular neighborhood of $J_i$, let $m(J_i)$ be a meridian of $J_i$ in $\partial T(J_i)$, and 
let $l(J_i)$ be a parallel of $J_i$ in $\partial T(J_i)$. 
Define the manifold $M=\chi (M_0,L)$ obtained from $M_0$ by surgery along the framed link $L$
by $M=X_L\amalg_h (\amalg_{i\in I} T_i)$, where $X_L=M_0\setminus \amalg_{i\in I} T(J_i)$ is the exterior of $L$,
the $T_i$ are solid tori,
and $h:\partial X_L =\amalg_{i\in I} \partial T(J_i) \to \amalg_{i\in I} \partial T_i$ is a homeomorphism sending 
$\partial T(J_i)$ onto $\partial T_i$ and the parallel $l(J_i)$ onto a meridian $m_i$ of $T_i$, for each $i\in I$.
We have :
$$H_1(\chi(M_0,L),\mathbb{Q})=\frac{\bigoplus_{i\in I} \mathbb{Q}\,m(J_i)}{\bigoplus_{i\in I} \mathbb{Q}\,l(J_i)}$$
\paragraph{Remark} If $H_1(M_0,\mathbb{Z})=0$, $\mathbb{Q}$ can be replaced by $\mathbb{Z}$.

In the manifold $M$, one can consider the link $\hat{L}$, whose components are the cores of the tori $T_i$, parallelised by 
the meridians $m(J_i)$. The surgery along $\hat{L}$ in $M$ is the \emph{inverse surgery} of the surgery along $L$ in $M_0$.
We have $\chi(\chi(M_0,L),\hat{L})=M_0$.

A \emph{surgery presentation} is a triple
$(M_0,L,K_0)$, where $M_0$ is a $\mathbb{Q}$-sphere, $K_0$ is a trivial knot in $M_0$, and $L=\amalg_{i=1}^n J_i$ is a 
framed link in $M_0$, such that $lk(K_0,J_i)=0$ for $1\leq i\leq n$. 
It is called \emph{admissible} if $\det((lk(J_i,J_j))_{1\leq i,j\leq n})\neq 0$. 
\begin{lemma} Set $(M,K)=(\chi(M_0,L),K_0)$.
\begin{enumerate}
 \item The knot $K$ is null-homologous. 
 \item The manifold $M$ is a $\mathbb{Q}$-sphere if and only if the surgery presentation is admissible.
\end{enumerate}
\end{lemma}
\proof
Consider a disk $D$ bounded by $K_0$, and tranverse to all the $J_i$. For all $i$, $<J_i,D>=0$. 
So $K$ bounds a surface obtained from $D$ by tubing, that is by adding tubes around the $J_i$ connecting pairs of points 
of $J_i\cap D$ with opposite signs.

We have :
$$H_1(M;\mathbb{Q})=\frac{\bigoplus_{i=1}^n \mathbb{Q}\,m(J_i)}{\bigoplus_{i=1}^n \mathbb{Q}\,l(J_i)}.$$
Now $l(J_i)=\sum_{j=1}^n lk(J_i,J_j) m(J_j)$, 
so $(lk(J_j,J_i))_{1\leq i,j\leq n}$ is a presentation matrix for the $\mathbb{Q}$-module $H_1(M;\mathbb{Q})$.
Thus $\det((lk(J_j,J_i))_{1\leq i,j\leq n})\neq 0$ if and only if $H_1(M;\mathbb{Q})=0$.\fin

Such a presentation always exists :
\begin{lemma}
 If $K$ is a null-homologous knot in a $\mathbb{Q}$-sphere $M$, then $K$ admits a surgery presentation $(S^3,L,K_0)$.
\end{lemma}
See Garoufalidis and Kricker \cite[p.117, Fact 1]{G-K} for the case $M=S^3$.
The generalization to any $\mathbb{Q}$-sphere $M$ easily follows from :
\begin{theorem}[Lickorish-Wallace,1960] \label{LW}
Any 3-manifold without boundary can be obtained from $S^3$ by surgery along a framed link.
\end{theorem}

  \subsection{Hermitian presentation matrix} \label{sechpm}

\begin{definition}
 A matrix $A\in \mathcal{M}(\qtt)$ is called \emph{hermitian} if $\bar{A}^t=A$, where $\bar{A}$
is defined by $\bar{A}(t)=A(t^{-1})$.
\end{definition}

Given a null-homologous knot $K$ in a rational homology sphere $M$, we aim to determine a hermitian
presentation matrix of the Alexander module $\mathcal{A}(K)$.
\begin{proposition} \label{propmatpres}
 Let $K$ be a null-homologous knot in a rational homology sphere $M$. Consider a surgery presentation
$(M_0,L=\amalg_{i=1}^n J_i,K_0)$ of $K$. Let $\tilde{X}_0$ be the infinite cyclic covering 
of the exterior $X_0$ of $K_0$, and, for $1\leq i\leq n$, let $\tilde{J}_i$ be a lift of $J_i$ in $\tilde{X}_0$. 
The equivariant linking matrix $A(t)$ defined by $A_{ij}(t)=lk_e(\tilde{J}_j,\tilde{J}_i)$ is a presentation
matrix of the Alexander module $\mathcal{A}(K)$. This matrix is hermitian and satisfies $\det(A(1))\neq 0$.
\end{proposition}

\proof
We will need the following lemma :
\begin{lemma} \label{lemmaHi}
 Let $U$ be a trivial knot in a $\mathbb{Q}$-sphere $M$, let $X$ be the exterior of $U$, and let $\tilde{X}$ 
 be the infinite cyclic covering associated. We have $H_1(\tilde{X},\mathbb{Q})=0$ and $H_2(\tilde{X},\mathbb{Q})=0$.
\end{lemma}
\proof
Let $D$ be a disk bounded by $U$, and set $Y=M\setminus D$. Denote by $\tau$ the automorphism of $\tilde{X}$
corresponding to the action of $t$. Let $\tilde{Y}$ be (resp. $\tilde{D}$) a copy of $Y$ (resp. $D$) in $\tilde{X}$. 
Set $Y_p=\amalg_{k\in\mathbb{Z}} \tau^{2k}(Y)$, and $Y_i=\amalg_{k\in\mathbb{Z}} \tau^{2k+1}(Y)$.

We first show that $H_2(\tilde{X},\mathbb{Q})=0$.
Use the Mayer-Vietoris sequence associated with $\tilde{X}=Y_p\cup Y_i$ :
$$H_2(\amalg_{k\in\mathbb{Z}} \tau^k(\tilde{D})) \longrightarrow H_2(Y_p)\oplus H_2(Y_i) \longrightarrow 
H_2(\tilde{X}) \longrightarrow H_1(\amalg_{k\in\mathbb{Z}} \tau^k(\tilde{D}))$$
For $i>0$, $H_i(\amalg_{k\in\mathbb{Z}} \tau^k(\tilde{D}))=\qtt\otimes_{\mathbb{Q}} H_i(D)=0$.
Thus $H_2(\tilde{X})\simeq H_2(Y_p)\oplus H_2(Y_i) \simeq \qtt \otimes_{\mathbb{Q}} H_2(Y)$. 
It remains to show that $H_2(Y)=0$.

Let $V$ be a regular neighborhood of $D$. Since $Y$ deformation retracts to $Z=M\setminus V$, 
we have $H_2(Y)=H_2(Z)$.
Then use the exact sequence associated with $(M,Z)$.
$$H_3(M)=\mathbb{Q}[M] \longrightarrow H_3(M,Z) \longrightarrow H_2(Z) \longrightarrow H_2(M)=0$$
By excision, $H_3(M,Z)=H_3(V,\partial V)=\mathbb{Q}[V]$. Since $[V]$ is the image of $[M]$ by the first morphism,
the second morphism is trivial and $H_2(Z)=0$.

We now show, with the same method, that $H_1(\tilde{X},\mathbb{Q})=0$. 
Use the Mayer-Vietoris sequence associated with $\tilde{X}=Y_p\cup Y_i$ :
$$0 \longrightarrow H_1(Y_p)\oplus H_1(Y_i) \longrightarrow H_1(\tilde{X}) \longrightarrow 
H_0(\amalg_{k\in\mathbb{Z}} \tau^k(\tilde{D})) \longrightarrow H_0(Y_p)\oplus H_0(Y_i)$$
The last morphism is injective, so
$H_1(\tilde{X})\simeq H_1(Y_p)\oplus H_1(Y_i) \simeq \qtt \otimes_{\mathbb{Q}} H_1(Y)$. 
It remains to show that $H_1(Y)=0$. We have $H_1(Y)=H_1(Z)$. Then use the exact sequence associated with
the pair $(M,Z)$.
$$H_2(M)=0 \longrightarrow H_2(M,Z) \longrightarrow H_1(Z) \longrightarrow H_1(M)=0$$
Thus $H_1(Z) \simeq H_2(M,Z) \simeq H_2(V,\partial V)$. Since $V$ is a ball, we have $H_2(V,\partial V)=0$. \fin

We now prove the proposition.
First note that the $J_i$ lift in $\tilde{X}_0$ into homeomorphic copies, because $lk(J_i,K_0)=0$.
The infinite cyclic covering $\tilde{X}$ associated with $K\subset M$ is obtained from $\tilde{X}_0$ 
by surgery along all the $\tilde{J}_{i,k}=\tau_0^k(\tilde{J}_i)$, where $\tau_0$ is the generator of $Aut(\tilde{X}_0)$
which induces the action of $t$. 

We calculate $H_1(\tilde{X};\mathbb{Q})$,
using the exact sequence associated with the pair $(\tilde{X},\tilde{Z})$, where $\tilde{Z}$ is the preimage 
of $Z=M_0\setminus \amalg_{i=1}^n Int(T(J_i))$ by the covering map. 
$$H_2(\tilde{X},\tilde{Z}) \fl{\partial}  H_1(\tilde{Z}) \longrightarrow H_1(\tilde{X})
\longrightarrow H_1(\tilde{X},\tilde{Z})$$
By excision, we have $H_j(\tilde{X},\tilde{Z})=H_j(\amalg_{1\leq i\leq n,k\in\mathbb{Z}} \tilde{T}_{i,k},\amalg_{1\leq i\leq n,k\in\mathbb{Z}}\partial\tilde{T}_{i,k})$ 
for all $j$, where the $\tilde{T}_{i,k}$ are the lifts of the solid tori $T_i$ glued during
the surgery. Thus $H_j(\tilde{X},\tilde{Z})=\bigoplus_{1\leq i\leq n,k\in\mathbb{Z}} 
H_j(\tilde{T}_{i,k},\partial \tilde{T}_{i,k})$, and we have $H_1(\tilde{X},\tilde{Z})=0$.
For $j=2$, we get $H_2(\tilde{X},\tilde{Z})=\bigoplus_{1\leq i\leq n,k\in\mathbb{Z}} \mathbb{Q}\ D_{i,k}=
\bigoplus_{1\leq i\leq n} \qtt\ D_{i,0}$, where $D_{i,k}$ is a meridian disk for 
$\tilde{T}_{i,k}$ bounded by $l(\tilde{J}_{i,k})$.

We now calculate $H_1(\tilde{Z};\mathbb{Q})$ using the exact sequence associated with the pair $(\tilde{X}_0,\tilde{Z})$.
$$H_2(\tilde{X}_0) \longrightarrow H_2(\tilde{X}_0,\tilde{Z}) \longrightarrow H_1(\tilde{Z}) 
\longrightarrow H_1(\tilde{X}_0)$$
By Lemma \ref{lemmaHi}, $H_2(\tilde{X}_0,\tilde{Z}) \rightarrow H_1(\tilde{Z})$ is an isomorphism.
Moreover, $$ H_2(\tilde{X}_0,\tilde{Z}) = \bigoplus_{1\leq i\leq n,k\in\mathbb{Z}} \mathbb{Q}\ \Sigma_{i,k},$$ 
where $\Sigma_{i,k}$ is a meridian disk for $T(\tilde{J}_{i,k})$ bounded by $m(\tilde{J}_{i,k})$. So :
$$H_1(\tilde{Z})=\bigoplus_{1\leq i\leq n,k\in\mathbb{Z}} \mathbb{Q}\ m(\tilde{J}_{i,k}) = 
\bigoplus_{1\leq i\leq n} \mathbb{Q} [t,t^{-1}]\ m(\tilde{J}_i).$$
We can then rewrite the first sequence :
$$\bigoplus_{1\leq i\leq n} \qtt\ D_{i,0} \fl{\partial} \bigoplus_{1\leq i\leq n} 
\mathbb{Q} [t,t^{-1}]\ m(\tilde{J}_i) \longrightarrow H_1(\tilde{X};\mathbb{Q}) \longrightarrow 0$$
We get a presentation of $H_1(\tilde{X};\mathbb{Q})$ with $n$ generators $m(\tilde{J}_i)$ and $n$ relations 
given by the images of $D_{i,0}$, that are the parallels $l(\tilde{J}_i)$. 

We now write $l(\tilde{J}_i)$ in function of the $m(\tilde{J}_i)$.
$$\begin{array}{lll}
l(\tilde{J}_i) &=& \displaystyle\sum_{j=1}^n \sum_{k\in\mathbb{Z}} lk(\tilde{J}_i,\tilde{J}_{j,k})\ m(\tilde{J}_{j,k}) \\
	&=& \displaystyle\sum_{j=1}^n \sum_{k\in\mathbb{Z}} lk(\tilde{J}_i,\tilde{J}_{j,k})\ t^k\ m(\tilde{J}_j) \\
	&=& \displaystyle\sum_{j=1}^n lk_e(\tilde{J}_i,\tilde{J}_j)\ m(\tilde{J}_j)
\end{array} $$
Therefore $A(t)$ is a presentation matrix for $\mathcal{A}(K)$.

It follows from the properties of the equivariant linking number that $A(t)$ is hermitian.
Since $M$ is a $\mathbb{Q}$-sphere, we have $\det(A(1))=\det((lk(J_j,J_i))_{1\leq i,j\leq n})\neq 0$. \fin

  \subsection{Properties of the Alexander module}

Consider a null-homologous knot $K$ in a $\mathbb{Q}$-sphere $M$, and write its Alexander module 
$\mathcal{A}(K)=\bigoplus_{i=1}^p \frac{\qtt}{\delta_i(t)}$,
with $\delta_{i+1}|\delta_i$ for $1\leq i<p$. Denote by $\doteq$ the equality modulo a unit of 
$\qtt$.
\begin{lemma}
For $1\leq i \leq p$, $\delta_i(t)\doteq \delta_i(t^{-1})$ and $\delta_i(1)\neq 0$.
\end{lemma}
\proof
Consider an equivariant linking matrix $A(t)$ associated with $K$.
By Proposition \ref{propmatpres}, $A(t)$ is a presentation matrix of $\mathcal{A}(K)$.
The matrix $A$ is equivalent to a diagonal matrix $D$ with diagonal 
($1, \dots, 1, \delta_p, \dots, \delta_1$) : there are two matrices $P$ and $Q$ in 
$GL_n(\qtt)$ such that $A=PDQ$. So $A=\bar{A}^t=\bar{Q}^t \bar{D} \bar{P}^t$, and $A$ is also
equivalent to $\bar{D}$. By unicity of the $\delta_i$ modulo a unit of $\qtt$, we have 
$\delta_i(t)\doteq\delta_i(t^{-1})$.

Moreover, $\det(A(t))=\prod_{i=1}^p \delta_i(t)$, and this determinant does not vanish for $t=1$.~\fin

We call {\em degree\/} of a polynomial $\sum_{i=q}^r \alpha_i t^i \in \mathbb{Q} [t,t^{-1}]$ with $\alpha_q\neq 0$
and $\alpha_r\neq 0$ the integer $r-q$. The $\delta_i$ with even degree can be normalised as symmetric polynomials.
\begin{lemma}
For $1\leq i\leq p$, if $\delta_i$ has even degree $2m$, then there are $r_j\in\mathbb{Q}$ for $0\leq j\leq m$ such that :
$$\delta_i(t)\doteq r_0 + \sum_{j=1}^m r_j (t^j+t^{-j})$$
\end{lemma}
\proof
Multiplying by a suitable power of $t$, we get $\delta_i(t)\doteq \sum_{j=-m}^m r_j t^j$,
with $r_mr_{-m}\neq 0$. Since $\delta_i(t)\doteq \delta_i(t^{-1})$, there is $\lambda\in\mathbb{Q}$
and $k\in\mathbb{Z}$ such that $\sum_{j=-m}^m r_j t^j=\lambda t^k \sum_{j=-m}^m r_j t^{-j}$. 
Necessarily $k=0$, and $r_j=\lambda r_{-j}$ for $-m\leq j\leq m$. Thus $r_m=\lambda r_{-m} =\lambda^2 r_m$,
and $\lambda^2=1$. If $\lambda=-1$, then $\delta_i(1)=\sum_{j=-m}^m r_j=r_0=0$. But $\delta_i(1)\neq 0$, so
$\lambda=1$. \fin

The degree parity of $\delta_i$ is related to the multiplicity of the root $-1$. Indeed, since 
$\delta_i(t)\doteq \delta_i(t^{-1})$, if $\alpha$ is a root for $\delta_i$, then $\alpha^{-1}$
is a root for $\delta_i$ with the same multiplicity.
So the roots of $\delta_i$ come by pairs $(\alpha,\alpha^{-1})$, except for $\alpha=\pm1$.
The case $\alpha=1$ does not occur, so $\delta_i$ has odd degree if and only if 
it admits $-1$ as a root with odd multiplicity.

\begin{proposition} \label{propptes}
Consider a null-homologous knot $K$ in a $\mathbb{Q}$-sphere $M$, and its Alexander module 
$\mathcal{A}(K)=\bigoplus_{i=1}^p \frac{\qtt}{\delta_i(t)}$,
with $\delta_{i+1}|\delta_i$ for $1\leq i<p$. The $\delta_i$ satisfy the following conditions :
\begin{itemize}
	\item $\delta_i(t^{-1})\doteq \delta_i(t)$ and $\delta_i(1)\neq 0$ for $1\leq i\leq p$,
	\item if, for $1\leq i\leq p$, $m_i$ denotes the multiplicity of $-1$ as a root of $\delta_i$,
		then for any odd integer $m$, the number of indices $i$ such that $m_i=m$ is even.
\end{itemize}
\end{proposition}
\paragraph{Remark}
These conditions imply that the Alexander polynomial $\prod_{i=1}^p \delta_i$ has even degree.

\proof
Proposition \ref{propmatpres} gives a hermitian presentation matrix $A\in\mathcal{M}_n(\qtt)$ 
for $\mathcal{A}(K)$.
There are matrices $P,Q\in GL_n(\qtt)$ such that $D=PAQ$ is a diagonal matrix with
diagonal $(1, \dots, 1, \delta_p, \dots, \delta_1)$. The matrix 
$B=\bar{Q}^t A Q = \bar{Q}^t P^{-1} D$ is a hermitian matrix which has its $i^{th}$ column (and thus its 
$i^{th}$ row) divisible by $\delta_i$. Moreover, $\det(B)\doteq \det(A)\doteq\prod_{i=1}^p \delta_i$.
The multiplicity of $-1$ as a root for $\det(B)$ is $m_0=\sum_{i=1}^p m_i$. Denote by $\beta_{ij}$ the coefficients of $B$.
Note that $(t+1)^{max(m_i,m_j)}$ divides $\beta_{ij}$. 
\begin{eqnarray*}
\det(B) &=& \sum_{\sigma}sign(\sigma)\prod_{i=1}^n\beta_{i\sigma(i)} \\
       &=& \sum_{\sigma; m_i=m_{\sigma(i)} \forall i}sign(\sigma)\prod_{i=1}^n\beta_{i\sigma(i)}\qquad mod\ (t+1)^{m_0+1}
\end{eqnarray*}
Set $I_m=\{i\in\{1\dots n\} | m_i=m\}$, $J=\{ m\in\mathbb{N} | I_m\neq \emptyset \}$, and $B_m=(\beta_{ij})_{(i,j)\in I_m\times I_m}$.
\begin{eqnarray*}
\det(B) &=& \prod_{m\in J} \det(B_m) \qquad mod\ (t+1)^{m_0+1}
\end{eqnarray*}
For $m\in J$, the matrix $B_m$ is hermitian, and has all
coefficients divisible by $(t+1)^m$. So $\det(B_m)$ is symmetric and divisible by $(t+1)^{mk(m)}$, where $k(m)=|I_m|$. 
As $\det(B_m)$ is symmetric, it is divisible by an even power of $(t+1)$, so if
$mk(m)$ is odd, $\det(B_m)$ is divisible by $(t+1)^{mk(m)+1}$. But $\sum_{m\in J} mk(m) = m_0$,
and $\det(B)\neq 0\ mod\ (t+1)^{m_0+1}$, so for $m\in J$, the multiplicity of $-1$ as a root of
$\det(B_m)$ is $mk(m)$, and $mk(m)$ is even. Thus, if $m$ is odd, $k(m)$ is even. \fin
\paragraph{Remark}
This proposition is a consequence of the existence of the Blanchfield form (see Theorem \ref{thdecomposition}).
It shows one implication of Theorem \ref{Amodules}.

  \subsection{Realization of equivariant linking matrices} \label{secreal}

In this subsection, we prove :
\begin{proposition} \label{propreal}
Consider a matrix $A(t)$ with coefficients in $\qtt$, hermitian, and satisfying $\det(A(1))\neq 0$.
The matrix $A(t)$ is the equivariant linking matrix associated with a surgery presentation of a 
null-homologous knot $K$ in a $\mathbb{Q}$-sphere $M$.
\end{proposition}
\proof
We first need to realize arbitrary linking numbers in $\mathbb{Q}$-spheres.
\begin{lemma} \label{Lens}
 Consider coprime integers $n$ and $k$, with $n>0$. In the lens space $L(n,k)$ obtained by $(-\frac{n}{k})$-surgery 
on the trivial knot $U$, we have $H_1(L(n,k);\mathbb{Z})=\frac{\mathbb{Z}}{n\mathbb{Z}} m(U)$ and $lk(m(U),m(U))=\frac{k}{n}$.
\end{lemma}
\proof
Consider a tubular neighborhood $T(U)$ of $U$, and the exterior $X$ of $U$.
Define the preferred parallel $l_0(U)$ of $U$ as the intersection of $T(U)$ with a disk bounded by $U$, 
and define the surgery curve of $U$ by $l(U)=n\,m(U)-k\,l_0(U)$ in $H_1(\partial T(U))$. Now consider
the manifold $L(n,k)=\chi(S^3,U)$.
We have $H_1(L(n,k);\mathbb{Z})=\frac{\mathbb{Z}\,m(U)}{\mathbb{Z}\,l(U)}$. In $H_1(X)$, 
$l(U)=n\,m(U)$. Thus $H_1(L(n,k);\mathbb{Z})=\frac{\mathbb{Z}}{n\mathbb{Z}} m(U)$.

We now compute $lk(m(U),m(U))$.
Let $T=S^1\times D^2$ be the solid torus glued during the surgery. Set $S_t=\{x\in D^2 |\ ||x||=t\}$.
Define a homeomorphism $h : \partial T = S^1\times S_1 \rightarrow S^1\times S_{\frac{1}{2}}$
by $(s,x) \mapsto (s,\frac{1}{2} x)$, and set $m_{int}(U)=h(m(U))$. Then 
$lk(m(U),m(U))=lk(m_{int}(U),m(U))$. Consider a meridian disk $D$ of $T$ bounded by $l(U)$. 
Then $$lk(m_{int}(U),l(U))=<m_{int}(U),D>=<l(U),m(U)>_{\partial T(U)} =k,$$
and 
\begin{eqnarray*}
lk(m_{int}(U),l(U)) &=& lk(m_{int}(U),n\,m(U)-k\,l_0(U)) \\
  &=& n\,lk(m_{int}(U),m(U))-k\,lk(m_{int}(U),l_0(U)). 
\end{eqnarray*}
Since $lk(m_{int}(U),l_0(U))=0$, we have $lk(m(U),m(U))=\frac{k}{n}$. \fin

\begin{corollary} \label{step 1}
 Take $m>0$, and for all ${1\leq i \leq m}$, consider coprime integers $k_i\in \mathbb{Z}$ and $n_i\in \mathbb{N}^*$. 
There exist a $\mathbb{Q}$-sphere $M$ and curves $c_i$
such that $H_1(M;\mathbb{Z})=\oplus_{i=1}^m \mathbb{Z}/n_i\mathbb{Z}\ c_i$,
with $lk(c_i,c_i)=\frac{k_i}{n_i}\ mod\ \mathbb{Z}$ et $lk(c_i,c_j)=0$ pour $i\neq j$.
\end{corollary}
\proof
By Lemma \ref{Lens}, there are $m$ $\mathbb{Q}$-spheres $M_i$ such that
$H_1(M_i;\mathbb{Z})=\mathbb{Z}/n_i\mathbb{Z}\,c_i$, with $lk(c_i,c_i)=\frac{k_i}{n_i}$.
Define $M$ as the connected sum of the $M_i$. \fin

\begin{corollary} \label{step 2}
For any family of rational numbers $(a_{ij})_{1\leq j\leq i\leq r}$, there are a $\mathbb{Q}$-sphere $M$ 
and simple, closed, pairwise disjoint, framed curves $f_i$, $1\leq i\leq r$, in $M$ such that 
$lk(f_i,f_j)=a_{ij}$ for $j\leq i$.
\end{corollary}
\proof
Given integers $(n_{ij})_{1\leq j\leq i \leq r}$ and $(k_{ij})_{1\leq j\leq i \leq r}$,
with $n_{ij}>0$ et $k_{ij}$ prime to $n_{ij}$ for all $1\leq j\leq i\leq r$,
Corollary \ref{step 1} gives a $\mathbb{Q}$-sphere $M$ such that 
$H_1(M;\mathbb{Z})=\bigoplus_{1\leq j\leq i \leq r} \mathbb{Z}/n_{ij}\mathbb{Z}\ c_{ij}$, with 
$lk(c_{ij},c_{ij})=\frac{k_{ij}}{n_{ij}}\ mod\ \mathbb{Z}$ and $lk(c_{ij},c_{st})=0$ if $(i,j)\neq (s,t)$. 
For $j<i$, take $n_{ij}$ and $k_{ij}$ such that $a_{ij}=\frac{k_{ij}}{n_{ij}}\ mod\ \mathbb{Z}$. 
For $i<j$, set $c_{ij}=c_{ji}$, $n_{ij}=n_{ji}$, and $k_{ij}=k_{ji}$.
Then, for $1\leq i \leq r$, choose curves $\gamma_i$
such that $\gamma_i=\sum_{l\neq i} c_{il}$ in $H_1(M;\mathbb{Z})$. For $i\neq j$, we get :
\begin{eqnarray*}
 lk(\gamma_i,\gamma_j) &=& lk(c_{ij},c_{ij})\ mod\ \mathbb{Z} \\
  &=& \frac{k_{ij}}{n_{ij}}\ mod\ \mathbb{Z}
\end{eqnarray*}
For $1\leq i\leq r$, choose $n_{ii}$ and $k_{ii}$ such that 
$lk(\gamma_i,\gamma_i)+\frac{k_{ii}}{n_{ii}}=a_{ii}\ mod\ \mathbb{Z}$, then consider $f'_i$
such that $f'_i=\gamma_i+c_{ii}$ in $H_1(M;\mathbb{Z})$. Then :
\begin{eqnarray*}
 lk(f'_i,f'_i) &=& lk(\gamma_i,\gamma_i)+lk(c_{ii},c_{ii})\ mod\ \mathbb{Z} \\
  &=& a_{ii}\ mod\ \mathbb{Z},
\end{eqnarray*}
and, if $j<i$, $$lk(f'_i,f'_j)=lk(\gamma_i,\gamma_j)=a_{ij}\ mod\ \mathbb{Z}.$$ 
Now, we can choose the curves $f_i$ in the homology classes of the $f'_i$, and their preferred parallels $l(f_i)$,
such that $lk(f_i,f_j)=a_{ij}$ in $\mathbb{Q}$ for all $1\leq j\leq i\leq r$.\fin

We now prove the proposition, generalizing the method of Levine \cite{Lev} for knots in $S^3$.

Write $A(t)=(P_{ji}(t))_{1\leq i,j\leq n}$, and $P_{ij}(t)=\sum_{k=-d}^d r_{ij}^{(k)} t^k$, with $r_{ij}^{(k)}=r_{ji}^{(-k)}$
for all $i,j,k$.
By Corollary \ref{step 2}, there are a $\mathbb{Q}$-sphere $M_0$ and pairwise disjoint simple closed framed curves
$\gamma_{ik}$, $1\leq i \leq n$, $0\leq k \leq d$, such that, for $1\leq i,j\leq n$ and $0\leq k,l\leq d$ :
$$lk(\gamma_{ik},\gamma_{jl})=\left\{ 
	\begin{array}{l l} r_{ij}^{(k)} & if\ l=0 \\ r_{ji}^{(l)} & if\ k=0 \\ 0 & otherwise \end{array} \right.$$

Consider a disk $D$ in $M_0$, disjoint from all $\gamma_{ik}$. Set $K_0=\partial D$.
For $1\leq i\leq n$ and $1\leq k\leq d$, connect $\gamma_{i,k-1}$ to $\gamma_{ik}$ with an arc $\alpha_{ik}$ 
such that $<\alpha_{ik},D>=1$, where $<.,.>$ denotes the algebraic intersection number.
Choose the $\alpha_{ik}$ pairwise disjoint, and disjoint from all the $\gamma_{ik}$.
Now consider a band $B_{ik}=h_{ik}([-1,1]\times [0,1])$ such that :
	\begin{itemize}
		\item $h_{ik}(\{0\}\times [0,1])=\alpha_{ik}$,
		\item $h_{ik}([-1,1]\times \{0\})=B_{ik}\cap (-\gamma_{i,k-1})\subset \partial B_{ik}$,
		\item $h_{ik}([-1,1]\times \{1\})=B_{ik}\cap \gamma_{ik}\subset -\partial B_{ik}$.
	\end{itemize}
We also suppose that $B_{ik}$ meets the tubular neighborhoods $T(\gamma_{i,k-1})$ and $T(\gamma_{ik})$ along
the preferred parallels $l(\gamma_{i,k-1})$ and $l(\gamma_{ik})$.

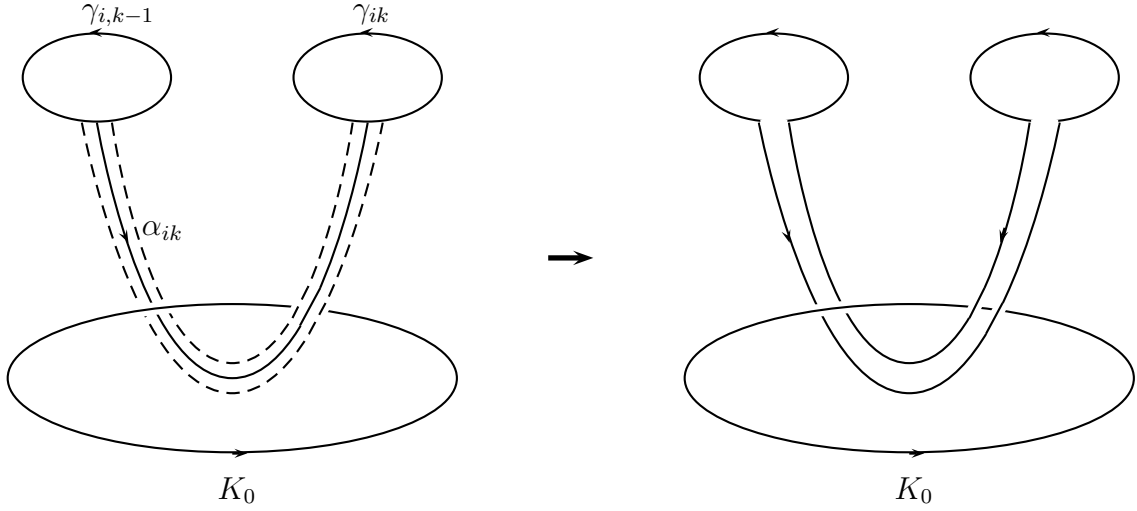
\begin{figure}[htb] 
\begin{center}
\begin{pspicture}(0,0.4)(15,7)
\psellipse(1.2,6)(1,0.6)
\psellipse(4.8,6)(1,0.6)
\pscurve(1.2,5.4)(3,2)(4.8,5.4)
\pscurve[linestyle=dashed](1.4,5.4)(3,2.2)(4.6,5.4)
\pscurve[linestyle=dashed](1,5.4)(3,1.8)(5,5.4)

\psellipse[border=2pt](3,2)(3,1)
\psline[border=2pt,linestyle=dashed](3.8,2.8)(4,3.2)
\psline[border=2pt](3.9,2.7)(4.17,3.2)
\psline[border=2pt,linestyle=dashed](4.02,2.59)(4.24,3)

\psline[linewidth=1pt]{<-}(1.1,6.6)(1.2,6.6)
\put(1,6.8){$\gamma_{i,k-1}$}
\psline[linewidth=1pt]{<-}(4.7,6.6)(4.8,6.6)
\put(4.6,6.8){$\gamma_{ik}$}
\psline[linewidth=1pt]{->}(3,1)(3.2,1)
\put(2.8,0.4){$K_0$}
\psline[linewidth=1pt]{->}(1.57,3.9)(1.62,3.8)
\put(1.8,3.9){$\alpha_{ik}$}

\psline[linewidth=2pt]{->}(7.2,3.6)(7.8,3.6)

\psellipse(10.2,6)(1,0.6)
\psellipse(13.8,6)(1,0.6)
\pscurve(10.4,5.4)(12,2.2)(13.6,5.4)
\pscurve(10,5.4)(12,1.8)(14,5.4)

\psellipse[border=2pt](12,2)(3,1)
\psline[border=2pt](12.8,2.8)(13,3.2)
\psline[border=2pt](13.02,2.59)(13.24,3)
\psline[linecolor=white,linewidth=2pt](10.02,5.4)(10.38,5.4)
\psline[linecolor=white,linewidth=2pt](13.62,5.4)(13.98,5.4)

\psline[linewidth=1pt]{<-}(10.1,6.6)(10.2,6.6)
\psline[linewidth=1pt]{<-}(13.7,6.6)(13.8,6.6)
\psline[linewidth=1pt]{->}(12,1)(12.2,1)
\put(11.8,0.4){$K_0$}
\psline[linewidth=1pt]{->}(10.38,3.9)(10.41,3.8)
\psline[linewidth=1pt]{->}(13.27,4)(13.21,3.8)

\end{pspicture}
\end{center}
\end{figure}

For all $i$, define a knot $J_i$ by : 
	$$J_i=\lhp \lbp \coprod_{k=0}^d \gamma_{ik} \rbp \setminus \lbp \coprod_{k=0}^d h_{ik}([-1,1]\times \partial [0,1]) 
	\rbp \rhp \bigcup \lhp \coprod_{k=0}^d h_{ik}(\partial [-1,1] \times [0,1]) \rhp .$$
Define the preferred parallels $l(J_i)$ similarly from the $l(\gamma_{ik})$.
Note that $J_i$ is homologous to $\sum_{k=0}^d \gamma_{ik}$.

Set $L=\amalg_{i=1}^n J_i$. For all $i$, we have $lk(J_i,K_0)=0$. Since $\det(A(1))\neq 0$, $(M_0,L,K_0)$ is an admissible surgery presentation.
Set $(M,K)=(\chi(M_0,L),K_0)$.

The $J_i$ lift in the infinite cyclic covering $\tilde{X}_0$ associated with $K_0$. Since $<J_i,D>=0$, the lifts 
are homeomorphic to the $J_i$. Choose lifts $\tilde{J}_i$ such that the corresponding lifts of 
the $\gamma_{i0}$ all lie in the same copy of $M\setminus D$ in $\tilde{X}_0$.
By Proposition \ref{propmatpres}, the equivariant linking matrix associated with the surgery presentation 
$(M_0,L,K_0)$ is a square matrix of order $n$, whose coefficients are 
$lk_e(\tilde{J}_i,\tilde{J}_j)=\sum_{k\in\mathbb{Z}} lk(\tilde{J_i},\tilde{J}_{j,k}) t^k$, with $\tilde{J}_{j,k}=\tau_0^k(\tilde{J}_j)$.
We have $$lk(\tilde{J_i},\tilde{J}_{j,k})=\left\{
	\begin{array}{l l} \displaystyle \sum_{l=\max(0,k)}^{\min(d,d+k)} lk(\gamma_{il},\gamma_{j,l-k}) & if\ -d\leq k\leq d \\
		0 & otherwise \end{array} \right.$$
Fix $k$ such that $-d\leq k\leq d$. 
$$lk(\tilde{J_i},\tilde{J}_{j,k})=\left\{
	\begin{array}{l l} r_{ij}^{(k)} & if\ k\geq 0 \\ r_{ji}^{(-k)} & if\ k\leq 0 \end{array} \right.$$
Since $r_{ji}^{(-k)}=r_{ij}^{(k)}$, we get $lk_e(\tilde{J}_i,\tilde{J}_j)=P_{ij}(t)$.
\fin

  \subsection{Proof of Theorem \ref{Amodules}} \label{secth1}

The modules satisfying the conditions of Theorem \ref{Amodules} can be written as a direct sum of terms 
$\frac{\qtt}{(P)}$, with $P$ symmetric, $P(1)\neq 0$, and terms $\frac{\qtt}{((1+t)^n)}\oplus\frac{\qtt}{((1+t)^n)}$,
with $n$ odd. Hence we have to realize these two kinds of modules, and the direct sums of Alexander modules.

\begin{lemma} \label{realpol}
If $P\in \qtt$ satisfies $P(t^{-1})=P(t)$ et $P(1)\neq 0$, 
then there exist a $\mathbb{Q}$-sphere $M$ and a null-homologous knot $K$ in $M$ such that 
the Alexander module of $K$ is $\mathcal{A}(K)=\qtt/(P)$.
\end{lemma}
\proof Apply Proposition \ref{propreal} and Proposition \ref{propmatpres} to the $1\times 1$ matrix $A(t)=(P(t))$. \fin

\begin{lemma}
For all integer $n>0$, there exist a $\mathbb{Q}$-sphere $M$ and a null-homologous knot $K$ in $M$
such that the Alexander module of $K$ is : 
$$\mathcal{A}(K)=\frac{\qtt}{(1+t)^n} \oplus \frac{\qtt}{(1+t)^n}$$
\end{lemma}
\proof Apply Proposition \ref{propreal} and Proposition \ref{propmatpres} to the matrix :
$$A(t)=\begin{pmatrix} 0 & (1+t)^n \\ (1+t^{-1})^n & 0 \end{pmatrix}.$$ \fin

Given two null-homologous knots $K_1$ and $K_2$ in $\mathbb{Q}$-spheres 
$M_1$ and $M_2$, we define the connected sum $K_1\sharp K_2$ of $K_1$ and $K_2$ in the connected sum $M_1\sharp M_2$ of
$M_1$ and $M_2$ in the following way. For $i=1,2$, remove from $M_i$ a ball $B_i$ which intersects $K_i$ 
along an arc $\alpha_i([0,1])$, trivial in the sense that the complement of a regular neighborhood of
$\alpha_i([0,1])$ in $B_i$ is a solid torus. Orient these arcs from $\alpha_i(0)$ to $\alpha_i(1)$.
Then glue $M_1\setminus Int(B_1)$ to $M_2\setminus Int(B_2)$ along
a homeomorphism $h : \partial B_1 \to -\partial B_2$
such that $h(\alpha_1(0))=\alpha_2(1)$ et
$h(\alpha_1(1))=\alpha_2(0)$. 

\begin{lemma} \label{propsomme}
If $K_1$ and $K_2$ are null-homologous knots in $\mathbb{Q}$-spheres 
$M_1$ and $M_2$ respectively, then their connected sum $K=K_1\sharp K_2$ is a null-homologous knot
in the $\mathbb{Q}$-sphere $M=M_1\sharp M_2$, and $\mathcal{A}(K)=\mathcal{A}(K_1)\oplus \mathcal{A}(K_2)$.
\end{lemma}
\proof
For $i=1,2$, consider a surgery presentation $(S^3_i,L_i,J_i)$ of $K_i$.
Setting $S^3_0=S^3_1\sharp S^3_2$, $J=J_1\sharp J_2$ and $L=L_1\amalg L_2$, we get a surgery presentation
$(S^3_0,L,J)$ of $K$.
Proposition \ref{propmatpres} gives the presentation matrix 
$\begin{pmatrix} A_1(t) & 0 \\ 0 & A_2(t) \end{pmatrix}$ for $\mathcal{A}(K)$, 
where, for $i=1,2$, $A_i(t)$ is a presentation matrix for $\mathcal{A}(K_i)$. \fin

\section{Blanchfield forms}

  \subsection{An expression of the Blanchfield form} \label{secB}

Consider a null-homologous knot $K$ in a $\mathbb{Q}$-sphere $M$, with a surgery presentation $(M_0,L=\amalg_{i=1}^nJ_i,K_0)$.
By Proposition \ref{propmatpres}, the associated equivariant linking matrix $A(t)$ is a presentation matrix 
of $\mathcal{A}(K)$. The generators of $\mathcal{A}(K)$ associated with this presentation are the meridians $m_i$
of the $\tilde{J}_i$ (fixed lifts of the $J_i$ in the infinite cyclic covering $\tilde{X}_0$ associated with $K_0$).
\begin{lemma} \label{lemmaBf}
Set $L=-A(t)^{-1}$. Then $lk_e(m_i,m_j)=L_{ji}$.
\end{lemma}
Note that $lk_e(m_i,m_i)$ is well defined, because the $m_i$ are framed by the $\partial T(\tilde{J}_i)$.
\proof The parallel $l(\tilde{J}_i)$ is rationally null-homologous, so it bounds a rational chain $\Sigma_i^{(1)}$ in $\tilde{X}_0$, 
which can be chosen to intersect the $T(\tilde{J}_{j,k})$ along meridian disks.
Removing from $\Sigma_i^{(1)}$ its intersections with the interiors of the $T(\tilde{J}_{j,k})$ for all $j\in\{ 1\dots n\}$, $k\in\mathbb{Z}$, 
we get a new rational chain $\Sigma_i^{(2)}$ such that :
$$\partial \Sigma_i^{(2)}=l(\tilde{J}_i)-\sum_{j=1}^n lk_e(\tilde{J_i},\tilde{J_j})m_j.$$
The chain $\Sigma_i^{(2)}$ can be viewed in the infinite cyclic covering $\tilde{X}$ associated with $K$.
Now, in $\tilde{X}$, $l(\tilde{J}_i)$ bounds a disk $D_{i,0}$ contained in the solid torus $\tilde{T}_{i,0}$ glued
during the surgery. Gluing $D_{i,0}$ and $-\Sigma_i^{(2)}$ along $l(\tilde{J}_i)$, we get a rational chain
$\Sigma_i$ in $\tilde{X}$ such that :
$$\partial \Sigma_i=\sum_{j=1}^n lk_e(\tilde{J}_i,\tilde{J}_j)m_j=\sum_{j=1}^n A_{ji}(t) m_j.$$
This gives :
\begin{eqnarray*}
 lk_e(\partial\Sigma_i,m_j)&=&\sum_{k=1}^n A_{ki}(t) lk_e(m_k,m_j).
\end{eqnarray*}
Now, defining the $(m_j)_{int}$ as in the proof of Lemma \ref{Lens}, we have :
\begin{eqnarray*}
 <(m_j)_{int},\Sigma_i>_e 
  &=& \sum_{k\in\mathbb{Z}} <(m_j)_{int},\tau^k(\Sigma_i)>t^k \\
  &=& -\delta_{ij},
\end{eqnarray*}
where $<.,.>_e$ denotes the \emph{equivariant algebraic intersection number}, given by definition by the first equality.
So $lk_e(m_j,\partial \Sigma_i)=-\delta_{ij}$.
Thus $\sum_{k=1}^n A_{ki}(t) lk_e(m_k,m_j)=-\delta_{ij}$.\fin
\begin{corollary} \label{corBf}
 For $1\leq i,j\leq n$, $\phi_K(m_i,m_j)=(-A(t)^{-1})_{ji}\ mod\ \qtt$.
\end{corollary}

  \subsection{Orthogonal decomposition}

In this section, we prove Theorem \ref{thdecomposition}. We begin with three useful lemmas. 
\begin{lemma}
 If $\gamma \in \mathcal{A}$ has order $P$, then there exists $\eta\in \mathcal{A}$ such that
$\phi(\gamma,\eta)=\frac{1}{P}$.
\end{lemma}
\proof If there exists a strict divisor $Q$ of $P$ such that $Im(\phi(\gamma,.))\in \frac{1}{Q}\qtt$,
then $Q\gamma \in ker(\phi)$, and so $Q\gamma=0$, which is a contradiction. 
Thus $PIm(\phi(\gamma,.))$ is an ideal of $\qtt / (P)$ which contains a unit. So $PIm(\phi(\gamma,.))=\qtt / (P)$.\fin

\begin{lemma} \label{decomposition}
 Let $\gamma_1$, $\gamma_2 \in \mathcal{A}$ have respective orders $\pi_1$ and $\pi_2$. 
If $\pi_1$ is prime to $\overline{\pi}_2$, then $\phi(\gamma_1,\gamma_2)=0$.
\end{lemma}
\proof We have $\pi_1\phi(\gamma_1,\gamma_2) \in \qtt$ and
$\overline{\pi}_2 \phi(\gamma_1,\gamma_2) \in \qtt$. As $\qtt$ is a principal ideal domain, 
there are $A$ and $B$ in $\qtt$ such that $A\pi_1+B\overline{\pi}_2=1$,
thus $\phi(\gamma_1,\gamma_2)=A\pi_1\phi(\gamma_1,\gamma_2)+B\overline{\pi}_2\phi(\gamma_1,\gamma_2) \in
\qtt$. \fin

\begin{lemma} \label{directsum}
 Let $F$ be a $\qtt$-submodule of $\mathcal{A}$, and denote by $F^\perp$ its orthogonal with respect to $\phi$.
We have $\mathcal{A}=F\oplus F^\perp$ if and only if $F\cap F^\perp =0$.
\end{lemma}
\proof Write $F=\bigoplus_{i=1}^p \frac{\qtt}{(P_i)} x_i$. For $y\in\mathcal{A}$,
$\phi(x_i,y)=\frac{P}{P_i}$, with $P$ defined modulo $P_i$. So $\codim_{\mathbb{Q}}x_i^\perp\leq \deg(P_i)$ for $1\leq i\leq p$.
Now $F^\perp=\cap_{i=1}^p x_i^\perp$, so $\codim_{\mathbb{Q}} F^\perp \leq \sum_{i=1}^p \deg(P_i)= \dim_{\mathbb{Q}} F$.
Thus $\dim_{\mathbb{Q}}F+\dim_{\mathbb{Q}}F^\perp \geq \dim_{\mathbb{Q}} \mathcal{A}$. \fin

We now prove Theorem \ref{thdecomposition}. We can write 
$\mathcal{A}=\bigoplus_{i\in I} \frac{\qtt}{(\pi_i^{n_i})}$ with $I$ finite and each
$\pi_i$ prime. Lemma \ref{decomposition} gives an orthogonal decomposition into terms 
$\bigoplus_{i\in I} \frac{\qtt}{\pi^{n_i}}$
with $\pi$ prime and symmetric, and terms $(\bigoplus_{i\in I} \frac{\qtt}{(\pi^{n_i})}) \oplus
(\bigoplus_{i\in J} \frac{\qtt}{(\overline{\pi}^{n_i})})$, with $\pi$ prime and non symmetric.
We have three cases to treat (the case $\pi=1+t$ being particular).

\paragraph{First case} : $\mathcal{A}=\bigoplus_{i\in I} \frac{\qtt}{\pi^{n_i}}$, 
$\pi$ prime and symmetric, and $I$ finite. Denote $n=max\{n_i;i\in I\}$.

First note that there exists $\gamma\in\mathcal{A}$ such that $\phi(\gamma,\gamma)=\frac{P}{\pi^n}$, $P$
prime to $\pi$. Indeed, if such a $\gamma$ does not exist, then, for $\eta_1$ of order $\pi^n$ 
and $\eta_2$ such that $\phi(\eta_1,\eta_2)=\frac{1}{\pi^n}$, 
we have $\phi(\eta_1+\eta_2,\eta_1+\eta_2)=\frac{Q}{\pi^n}$ with $Q=2$ mod $\pi$, which gives a contradiction.

Consider such a $\gamma$, and denote by $F$ the $\qtt$-submodule of $\mathcal{A}$ generated by $\gamma$.
If $x\in F\cap F^\perp$, then $x=\lambda \gamma$ for a $\lambda\in\qtt$, and $\phi(x,\gamma)=0$ implies 
$\lambda\in\pi^n\qtt$.
Hence $\mathcal{A}=F\oplus F^\perp$.

As $dim_\mathbb{Q}(F^\perp)<dim_\mathbb{Q}(\mathcal{A})$, we can conclude by induction.

\paragraph{Second case} : $\mathcal{A}=(\bigoplus_{i\in I} \frac{\qtt}{(\pi^{n_i})}) \oplus
(\bigoplus_{i\in J} \frac{\qtt}{(\overline{\pi}^{n_i})})$, with $\pi$ prime, non symmetric, $\pi(-1)\neq 0$, 
and $I$ and $J$ finite. The conditions on $\pi$ imply that $\pi$ is prime to $\overline{\pi}$. In particular,
in the direct sum above, the two main terms are submodules isotropic for $\phi$. 

Denote $n=max\{n_i;i\in I\cup J\}$.
Without loss of generality, we suppose that $n$ appears as a power of $\pi$. Consider $\gamma_1$ of order $\pi^n$,
and $\gamma_2$ such that $\phi(\gamma_1,\gamma_2)=\frac{1}{\pi^n}$. Since projecting $\gamma_2$ on
$\bigoplus_{i\in J} \frac{\qtt}{(\overline{\pi}^{n_i})}$ does not change $\phi(\gamma_1,\gamma_2)$,
we can suppose $\gamma_2$ has order $\overline{\pi}^k$, and $\phi(\gamma_1,\gamma_2)=\frac{1}{\pi^n}$ implies $k=n$.

Now denote by $G$ the $\qtt$-submodule of $\mathcal{A}$ generated by $\gamma_1$ and $\gamma_2$.
If $x\in G\cap G^\perp$, write $x=\lambda \gamma_1+\mu\gamma_2$ with $\lambda,\mu\in\qtt$.
Then $\phi(x,\gamma_1)=0$ implies $\mu\in\overline{\pi}^n\qtt$, and $\phi(x,\gamma_2)=0$ implies $\lambda\in\pi^n\qtt$.
So $x=0$. Hence $\mathcal{A}=G\oplus G^\perp$. Again $dim_\mathbb{Q}(G^\perp)<dim_\mathbb{Q}(\mathcal{A})$, 
and we conclude by induction.

\paragraph{Third case} : $\mathcal{A}=\bigoplus_{i\in I} \frac{\qtt}{((1+t)^{n_i})}$, with $I$ finite.
Denote $n=max\{n_i;i\in I\}$.

If $n$ is even, we can replace $(1+t)^n$ by $(t+2+t^{-1})^{n/2}$, and proceed exactly like in the
first case. Now suppose $n$ is odd.

First note that, for all $\gamma\in\mathcal{A}$, $\phi(\gamma,\gamma)=\frac{P}{(1+t)^k}$, with $k$ even, $P(-1)\neq 0$
or $P=0$. Indeed, if $\phi(\gamma,\gamma)=\frac{P}{(1+t)^k}$, we have $\frac{P}{(1+t)^k}=\frac{\overline{P}}{(1+t^{-1})^k}$, 
so $P=t^k\overline{P}$, and $k$ odd implies $P(-1)=0$.

Consider $\gamma_1$ of order $(1+t)^n$ and $\gamma_2$ such that $\phi(\gamma_1,\gamma_2)=\frac{1}{(1+t)^n}$.
Note that $\gamma_2$ also has order $(1+t)^n$. Denote by $H$ the $\qtt$-submodule of $\mathcal{A}$ 
generated by $\gamma_1$ and $\gamma_2$. For $x\in H\cap H^\perp$, write $x=\lambda\gamma_1+\mu\gamma_2$, 
with $\lambda,\mu\in\qtt$. If $\lambda$ or $\mu$ is in $(1+t)^n\qtt$, then $x=0$. Suppose $\lambda$ and $\mu$ are not 
in $(1+t)^n\qtt$. We have :
$$\left\{ \begin{array}{l}
                \lambda \phi(\gamma_1,\gamma_1)+\frac{\mu}{(1+t^{-1})^n}=0 \\
		\frac{\lambda}{(1+t)^n}+\mu\phi(\gamma_2,\gamma_2)=0 \\
               \end{array} \right. .$$
Define minimal integers $j$ and $k$ such that $(1+t)^j\phi(\gamma_1,\gamma_1)$ and $(1+t)^k\phi(\gamma_2,\gamma_2)$
are in $\qtt$, and denote by $m_\lambda$ (resp. $m_\mu$) the multiplicity of the root $-1$ in $\lambda$ (resp. $\mu$).
We get :
$$\left\{ \begin{array}{l}
           m_\lambda-j=m_\mu-n \\
	   m_\lambda-n=m_\mu-k \\
          \end{array} \right. .$$
This implies $n-j=k-n$. But $n-j\geq 0$ and $k-n\leq 0$, so $j=k=n$, which is a contradiction, because n is odd, whereas 
$j$ and $k$ are even. Hence $\mathcal{A}=H\oplus H^\perp$.

We have showed in particular that $\lambda\gamma_1+\mu\gamma_2=0$ implies $(1+t)^n$ divides $\lambda$ and $\mu$.
So $H=\frac{\qtt}{((1+t)^n)}\gamma_1\oplus\frac{\qtt}{((1+t)^n)}\gamma_2$.
The last point is to describe the form $\phi$ over this submodule.

The matrix of the form $\phi$ on $H$, with respect to the basis $(\gamma_1, \gamma_2)$, is
$$\begin{pmatrix} \frac{\alpha(t)}{(t+2+t^{-1})^k} & \frac{1}{(1+t)^n} \\ \frac{1}{(1+t^{-1})^n} & 
  \frac{\beta(t)}{(t+2+t^{-1})^l} \end{pmatrix},$$ with $k,l<\frac{n}{2}$, $\alpha$ and $\beta$ symmetric.
We want to show that there is a basis such that $\alpha$ and $\beta$
vanish. We shall first get $\alpha=0$. 
Define $\eta_1=\gamma_1+\frac{(-1)^k\alpha(-1)}{2}(1+t)^{n-2k}\gamma_2$. We get :
\begin{eqnarray*}
 \phi(\eta_1,\eta_1) &=& \frac{\alpha(t)+\frac{(-1)^k\alpha(-1)}{2}(t^{k-n}+t^{n-k})+
 \frac{\alpha(-1)^2}{4}(t+2+t^{-1})^{(n-k-l)}\beta(t)}{(t+2+t^{-1})^{k}}\\
  \end{eqnarray*}
Since $n-k-l>0$, the numerator vanishes for $t=-1$. So $\phi(\eta_1,\eta_1)$ can be written 
$\phi(\eta_1,\eta_1)=\frac{\nu(t)}{(t+2+t^{-1})^{k'}}$, where $k'<k$. Now choose $\eta_2\in H$ such that 
$\phi(\eta_1,\eta_2)=\frac{1}{(1+t)^n}$. Replacing the basis $(\gamma_1, \gamma_2)$ by the basis $(\eta_1, \eta_2)$
makes the integer $k$ decrease. Iterating this process, we get a basis, again denoted $(\gamma_1, \gamma_2)$, such that $\alpha=0$.

Now consider $\gamma=\gamma_2+a(t)\gamma_1$. We have :
$$\phi(\gamma,\gamma)=\frac{a(t)+\overline{a}(t)t^n+\beta(t)t^l(1+t)^{n-2l}}{(1+t)^n}.$$
Set $a(t)=-\frac{1}{2} \beta(t)t^l(1+t)^{n-2l}$. Then $\overline{a}(t)t^n=a(t)$, so
$\phi(\gamma,\gamma)=0$. Hence in the basis $(\gamma_1, \gamma)$, we get $\alpha=0$ and $\beta=0$.
This concludes the proof of Theorem \ref{thdecomposition}.
\linebreak

As a direct consequence of Theorem \ref{thdecomposition}, we get the following lemma.
\begin{lemma} \label{lemmaremark}
For an Alexander module 
$\frac{\qtt}{(\pi^n)}\oplus\frac{\qtt}{(\overline{\pi}^n)}$, with $\pi$ prime and non symmetric, and $n$ odd if $\pi=(1+t)$, 
there is a unique isomorphism class of Blanchfield forms.
In particular, the realization of the Alexander module gives a realization of the Blanchfield form. 
\end{lemma}

Now, for an Alexander module 
$\mathcal{A}(K)=\frac{\qtt}{(\pi^n)}$ with $\pi$ prime and symmetric or $\pi=t+2+t^{-1}$, there may be different 
isomorphism classes of Blanchfield forms. The realization of the module described in Section \ref{secth1} gives a generator
$\gamma$ of $\mathcal{A}(K)$ for which $\phi_K(\gamma,\gamma)=\frac{-1}{\pi^n}$. Any other generator $\eta$ can be 
written $\eta=R\gamma$, with $R$ prime to $\pi$, and we have $\phi_K(\eta,\eta)=\frac{-R\bar{R}}{\pi^n}$. 
Here two questions arise. First, can any non degenerate hermitian form over such a module be realized as a Blanchfield 
form ? Second, what are the classes of symmetric polynomials prime to $\pi$ modulo $\pi^n$ and all the $R\bar{R}$ ?
The purpose of the next section is to give a positive answer to the first one. The last section gives a partial answer 
to the second one, showing there are infinitely many such classes.

  \subsection{Realization of Blanchfield forms}

\begin{proposition} \label{realBf}
Consider a $\qtt$-module $\mathcal{A}=\frac{\qtt}{(\Delta)}$ with $\Delta$ symmetric and $\Delta(1)\neq 0$.
Consider a non degenerate hermitian form $\phi : \mathcal{A}\times\mathcal{A} \to \frac{\mathbb{Q}(t)}{\qtt}$.
There is a null-homologous knot $K$ in a $\mathbb{Q}$-sphere $M$ such that $(\mathcal{A}(K),\phi_K)$
is isomorphic to $(\mathcal{A},\phi)$.
\end{proposition}
This result, together with Theorem \ref{thdecomposition}, Lemma \ref{lemmaremark}, and Lemma \ref{propsomme}, 
will conclude the proof of Theorem \ref{characterization}.
\proof[Proof of Proposition \ref{realBf}] Let $\gamma$ be a generator of $\mathcal{A}$. We have $\phi(\gamma,\gamma)=\frac{P}{\Delta}$,
with $P$ prime to $\Delta$.
Suppose there is a hermitian matrix $A(t)$ such that $\det(A)=r\Delta$ for an $r\in\mathbb{Q}$, 
and the cofactor $(1,1)$ of $A$ is $-rP$. By Proposition \ref{propreal} and Proposition \ref{propmatpres}, 
the matrix $A(t)$ is a presentation matrix of a null-homologous knot $K$ in a $\mathbb{Q}$-sphere $M$.
By Corollary \ref{corBf}, the first generator $m_1$ of this presentation satisfies $\phi_K(m_1,m_1)=\frac{P}{\Delta}$.
Let us construct such a matrix $A(t)$ to conclude.

Set :
$$A(t)=\begin{pmatrix} \alpha_1 & 1 &&&& \\ 1 & \ddots & \ddots && 0 & \\ & \ddots & \ddots & \ddots & & \\
	&& \ddots & \ddots & 1 & \\ & 0 && 1 & \alpha_i & q \\ &&&& q & \beta \end{pmatrix},$$
with $\alpha_j,\beta \in\qtt$ and $q\in\mathbb{Q}$. For $1\leq j\leq i$, consider the following sub-matrix of $A(t)$~: 

$$A_j(t)=\begin{pmatrix} \alpha_j & 1 &&&& \\ 1 & \ddots & \ddots && 0 & \\ & \ddots & \ddots & \ddots & & \\
	&& \ddots & \ddots & 1 & \\ & 0 && 1 & \alpha_i & q \\ &&&& q & \beta \end{pmatrix}.$$
The determinants of these matrices satisfy the following relations :
$$\left\{ \begin{array}{l l}
	\det(A_j)=\alpha_j \det(A_{j+1})-\det(A_{j+2}) & for\ 1\leq j\leq i-2 \\
	\det(A_{i-1})=\alpha_{i-1}\det(A_i)-\beta & \\
	\det(A_i)=\alpha_i\beta - q^2
	\end{array} \right. .$$
	
Now set $R_1=\Delta$ and $R_2=-P$. Consider the successive euclidean divisions $R_j=Q_jR_{j+1}-R_{j+2}$, with 
$d^\circ R_{j+2}<d^\circ R_{j+1}$. Let $k$ be the integer such that $r:=R_{k+2}$ is a non-zero rational number.
We want to identify these equalities with the above relations. As $r$ may not be a square, multiply the equalities 
by $r$. This gives $(rR_j)=Q_j(rR_{j+1})-(rR_{j+2})$. In particular, $(rR_k)=Q_k(rR_{k+1})-r^2$.
Set $i=k$, $q=r$, $\beta=rR_{k+1}$, and, for $1\leq j\leq k$, $\alpha_j=Q_j$. Then $\det(A_j)=rR_j$ for all $j$, and we get the required 
matrix $A(t)$. \fin

\subsection{Isomorphism classes of Blanchfield forms}

\begin{proposition}
 Consider the $\qtt$-module $\mathcal{A}=\frac{\qtt}{(\Delta^n)}$, with $\Delta$ prime and symmetric, or $\Delta(t)=t+2+t^{-1}$, and $n>0$.
The set $\mathcal{B}_{\Delta^n}$ of isomorphism classes of non degenerate hermitian forms on $\mathcal{A}$ is infinite.
\end{proposition}
This result and Theorem \ref{thdecomposition} prove Theorem \ref{thBf}.
\proof A non degenerate hermitian form $\phi$ on $\mathcal{A}$ is determined by the value of $\phi(\gamma,\gamma)$
for a generator $\gamma$ of $\mathcal{A}$. Consider two such forms $\phi_1$ and $\phi_2$, with $\phi_1(\gamma,\gamma)=\frac{P}{\Delta^n}$
and $\phi_2(\gamma,\gamma)=\frac{Q}{\Delta^n}$, where $P$ and $Q$ are symmetric, prime to $\Delta$, and defined modulo $\Delta^n$.
The forms $\phi_1$ and $\phi_2$ are isomorphic if and only if there is $R\in\qtt$ such that 
$\phi_1(R\gamma,R\gamma)=\phi_2(\gamma,\gamma)$, $i.e.$ $R\bar{R}P=Q\ mod\ \Delta^n$.

Symmetric polynomials can be written as polynomials in the variable $s=t+t^{-1}$. Define $\Delta_s$ by
$\Delta_s(t+t^{-1})=\Delta(t)$, and set $D_n=\frac{\mathbb{Q}[s]}{(\Delta_s^n)}$ and $E_n=\frac{\qtt}{(\Delta^n)}$.
Now consider the application $\varphi_n : E_n^\star \to D_n^\star$ given by $\varphi_n(R)(t+t^{-1})=R\bar{R}(t)$.
There is a natural bijection between $\mathcal{B}_{\Delta^n}$ and 
$coker(\varphi_n)=D_n^\star/\varphi_n(E_n^\star)$. The next lemma allows us to restrict our study to the case $n=1$.

\begin{lemma}
 There is an isomorphism $coker(\varphi_n) \fl{\sim} coker(\varphi_1)$.
\end{lemma}

\proof
 Consider the following commutative diagram of groups, where the morphisms $\pi_n$, $\pi_1$, $p_{\scriptscriptstyle{E}}$, 
$p_{\scriptscriptstyle{D}}$ are the natural projections :
\begin{center}
 \begin{tabular}{c c c c c}
  $E_n^\star$ & $\fl{\varphi_n}$ & $D_n^\star$ & $\fl{\pi_n}$ & $coker(\varphi_n)$ \\ & & & & \\
  $\quad \downarrow p_{\scriptscriptstyle{E}}$ & & $\quad \downarrow p_{\scriptscriptstyle{D}}$ & & \\ & & & & \\
  $E_1^\star$ & $\fl{\varphi_1}$ & $D_1^\star$ & $\fl{\pi_1}$ & $coker(\varphi_1)$ \\
 \end{tabular}
\end{center}
Let us show the existence of a morphism $p_c : coker(\varphi_n) \to coker(\varphi_1)$
such that the diagram still commutes. 
We have to show that $\pi_n(P)=\pi_n(Q)$ implies $\pi_1\circ p_{\scriptscriptstyle{D}}(P)=\pi_1\circ p_{\scriptscriptstyle{D}}(Q)$, 
\emph{i.e.} $ker(\pi_n)\subset ker(\pi_1\circ p_{\scriptscriptstyle{D}})$.
$$\pi_1\circ p_{\scriptscriptstyle{D}}(ker(\pi_n))=\pi_1\circ p_{\scriptscriptstyle{D}}(\varphi_n(E_n^\star))=
\pi_1\circ\varphi_1\circ p_{\scriptscriptstyle{E}}(E_n^\star)=0$$
So $p_c$ is well defined, and it is obviously surjective. It remains to show that it is injective.

Consider $x\in ker(p_c)$. There is $y$ in $D_n^\star$ such that $x=\pi_n(y)$. Then $\pi_1\circ p_{\scriptscriptstyle{D}}(y)=1$,
so $p_{\scriptscriptstyle{D}}(y)=\varphi_1\circ p_{\scriptscriptstyle{E}}(z)=p_{\scriptscriptstyle{D}}\circ \varphi_n(z)$. 
Thus $p_{\scriptscriptstyle{D}}(y\varphi_n(z^{-1}))=1$. It remains to show that 
$P=y\varphi_n(z^{-1})\in Im(\varphi_n)$.

We can write 
$P=1+\sum_{i=1}^{n-1} P_i\Delta_s^i$ with $d^\circ P_i<d^\circ \Delta_s$ for $1\leq i\leq n-1$. We shall search $R\in E_n^\star$
such that $R\bar{R}(t)=P(t+t^{-1})$ of the form $R=1+\sum_{i=1}^{n-1} R_i\Delta^i$, with $R_i$ symmetric and 
$d^\circ R_i<d^\circ \Delta$ for $1\leq i\leq n-1$. We have $R\bar{R}=1+\sum_{i=1}^{n-1} Q_i\Delta^i$, 
where $Q_i$ is the sum of $2R_i$ and a term depending on the $R_j$ for $j<i$. So we can define the $R_i$ by induction 
in order to get $Q_i(t)=P_i(t+t^{-1})$. Hence $P\in\varphi_n(E_n^\star)$.\fin

We now suppose $n=1$, and denote $D_1$, $E_1$, $\varphi_1$, by $D$, $E$, $\varphi$.
We first treat the case $\Delta(t)=t+2+t^{-1}$.
\begin{lemma}
 If $\Delta(t)=t+2+t^{-1}$, then $$coker(\varphi)\cong\mathbb{Q}^\star/(\mathbb{Q}^\star)^2
\cong\mathbb{Z}/2\mathbb{Z}\times\bigoplus_{\text{prime integers}}\mathbb{Z}/2\mathbb{Z}.$$
\end{lemma}
\proof
We have $D=\frac{\mathbb{Q}[s]}{(s+2)}\simeq\mathbb{Q}$ and $E=\frac{\qtt}{(t+2+t^{-1})}$. The elements of $E^\star$
can be written $at+b$, with $a\neq b$ in $\mathbb{Q}$. Now $\varphi(at+b)=a^2+b^2+abs=(a-b)^2\ mod\ (s+2)$.
So $Im(\varphi)=(\mathbb{Q}^\star)^2$. 

Each element of $\mathbb{Q}^\star$ can be written uniquely 
$\varepsilon\prod_{p\ prime} p^{n(p)}$ with only a finite number of non-zero $n(p)$, and $\varepsilon=\pm 1$.
Thus each element in the quotient 
$\mathbb{Q}^\star/(\mathbb{Q}^\star)^2$ can be written uniquely $\varepsilon\prod_{p\ prime} p^{n(p)}$ with $n(p)=1$ 
for a finite number of primes $p$, and $n(p)=0$ for the others.\fin

We now suppose $\Delta$ prime. In the general case, we will use some material of algebraic number theory. We will not detail 
these notions, and we refer the reader to \cite{Samuel}.

Since $\Delta$ is prime, $\Delta_s$ is prime too, so $D$ is a number field, that is a finite extension of $\mathbb{Q}$.
We can write $E=\frac{D[t]}{(t^2-st+1)}$, so $E$ is a Galois extension of $D$ of degree 2. We denote by $\sigma$ the 
only non trivial element of its Galois group. Denote by $\alpha$ the image of $t$ in $\frac{D[t]}{(t^2-st+1)}$. It is a primitive
element in $E$, \textit{i.e.} $E=D[\alpha]$, and the morphism $\sigma$ is given by $\sigma_{|D}=id_D$ and 
$\sigma(\alpha)=\alpha^{-1}$. Note that $\sigma$ is induced by the conjugation $t\mapsto t^{-1}$ in $\qtt$.
Define the norm $N_{E/D} : E \to D$ by $N_{E/D}(x)=x\sigma(x)$. On $E^\star$, the norm $N_{E/D}$ coincides with the morphism $\varphi$.

Denote by $A_D$ the ring of integers of $D$. A \emph{fractional ideal} of $D$ is a $A_D$-module $I$ such that $I\subset D$
and $kI\subset A_D$ for a $k\in D^\star$. The \emph{product} of two fractional ideals $I$ and $J$ is made of all finite
sums of products $ij$ with $i\in I$ and $j\in J$. It is again a fractional ideal. With this multiplication, 
the non trivial fractional ideals of $D$ form a group denoted by $\mathcal{F}_D$. The principal fractional ideals 
of $D$ are the fractional ideals of type $kA_D$ for $k\in D^\star$. They form a sub-group $\mathcal{P}_D$ of $\mathcal{F}_D$.
The quotient group $C_D=\mathcal{F}_D/\mathcal{P}_D$ is called the \emph{ideal class group} of $D$. 
\begin{theorem}[Dirichlet] (\cite[p.71, Theorem 2]{Samuel})\\
The ideal class group of a number field is finite.
\end{theorem}
We have a short exact sequence :
$$ 0 \to \mathcal{P}_D \to \mathcal{F}_D \to \mathcal{C}_D \to 0 . $$
This can be related to $D^\star$ by the surjective map $D^\star \twoheadrightarrow \mathcal{P}_D$ given by $k\mapsto kA_D$.

The field $E$ also is a number field, so all the previous definitions and results apply to $E$.
We shall define commutative diagrams based on the above exact sequence and surjective map, and the norm morphism.
For an ideal $xA_E$, define $N_{E/D}(xA_E)=N_{E/D}(x)A_D$. Now we have the following commutative diagram,
where the vertical arrows are given by the norm $N_{E/D}$.
\begin{center}
 \begin{tabular}{c c c}
  $E^\star$ & $\twoheadrightarrow$ & $\mathcal{P}_E$ \\
  $\downarrow$ & & $\downarrow$ \\
  $D^\star$ & $\twoheadrightarrow$ & $\mathcal{P}_D$ \\
 \end{tabular}
\end{center}
This induces a surjective map $coker(N_{E/D} : E^\star \to D^\star) \twoheadrightarrow coker(N_{E/D} : \mathcal{P}_E \to \mathcal{P}_D)$.
Recall our aim is to show that $coker(N_{E/D} : E^\star \to D^\star)$ is infinite. 
It suffices to show that $coker(N_{E/D} : \mathcal{P}_E \to \mathcal{P}_D)$ is infinite. 
To see this, we will use the above exact sequence.

For a fractional ideal $J$ of $E$, define $N_{E/D}(J)$ as the ideal of $D$ generated over $A_D$ by all the $N_{E/D}(x)$ 
for $x\in J$. For a principal ideal, we recover the previous definition. We get the following commutative diagram,
where the vertical arrows are again given by the norm $N_{E/D}$.
\begin{center}
 \begin{tabular}{c c c c c c c c c}
  0 & $\longrightarrow$ & $\mathcal{P}_E$ & $\longrightarrow$ & $\mathcal{F}_E$ & $\longrightarrow$ & $\mathcal{C}_E$ & $\longrightarrow$ & 0 \\
   & & $\downarrow$ & & $\downarrow$ & & $\downarrow$ & & \\
  0 & $\longrightarrow$ & $\mathcal{P}_D$ & $\longrightarrow$ & $\mathcal{F}_D$ & $\longrightarrow$ & $\mathcal{C}_D$ & $\longrightarrow$ & 0 \\
 \end{tabular}
\end{center}
Since $\mathcal{C}_E$ and $\mathcal{C}_D$ are finite, it remains to show that $coker(N_{E/D} : \mathcal{F}_E \to \mathcal{F}_D)$
is infinite.

Any fractional ideal $I$ of $D$ has a unique factorisation $I=\prod_{i=1}^q \mathfrak{p}_i^{e_i}$ where $q\in\mathbb{N}^\star$ 
and for each $i$, $\mathfrak{p}_i$ is a prime ideal of $A_D$ and $e_i\in\mathbb{Z}$ (\cite[p.60, Theorem 3]{Samuel}). 
Now, for any $J\in\mathcal{F}_D$, 
we have $N_{E/D}(JA_E)=J^2$ (\cite[p.25, Corollary 1]{Lang}). So $I$ is in the image of $N_{E/D}$ if and only if, 
for each $i$ such that $e_i$ is odd, $\mathfrak{p}_i$ is in the image of $N_{E/D}$.

Note that there are infinitely many prime ideals in $A_D$. Indeed, if $p$ is a prime integer, there is a prime ideal $I_p$
of $A_D$ such that $I_p\cap \mathbb{Z}=p\mathbb{Z}$ (\cite[p.9, Proposition 9]{Lang}).

Now consider a prime ideal $\mathfrak{p}$ of $A_D$ and set $J=\mathfrak{p}A_E$. The primes that appear in the decomposition 
of $J$ are exactly those whose intersection with $A_D$ is $\mathfrak{p}$. We have three cases (\cite[\S 5.2]{Samuel}) :
\begin{itemize}
 \item $J=\mathfrak{P}^2$, with $\mathfrak{P}$ prime in $A_E$, and $N_{E/D}(\mathfrak{P})=\mathfrak{p}$,
 \item $J$ itself is prime, and $N_{E/D}(J)=\mathfrak{p}^2$,
 \item $J=\mathfrak{P}_1\mathfrak{P}_2$, with $\mathfrak{P}_1$ and $\mathfrak{P}_2$ primes, $\mathfrak{P}_2=\sigma(\mathfrak{P}_1)$,
  and $N_{E/D}(\mathfrak{P}_i)=\mathfrak{p}$ for $i=1,2$.
\end{itemize}
In the second case, the ideal $\mathfrak{p}$ is not in the image of $N_{E/D}$. It follows from the Tchebotarev density theorem
(\cite[p.169, Theorem 10]{Lang}) that the density of the subset of prime ideals $\mathfrak{p}$ of $A_D$ for which $J$ is prime 
is equal to $\frac{1}{2}$ in the set of all prime ideals. This concludes the proof.\fin

\end{document}